\documentclass[11pt,twoside,reqno]{amsart}
\usepackage{mathptmx,amsmath,amssymb,amsfonts,amsthm,
mathptmx,enumerate,color}
\usepackage{graphicx}

\usepackage{epstopdf}

\newtheorem{theorem}{Theorem}
\newtheorem{corollary}{Corollary}
\newtheorem{lemma}{Lemma}

\theoremstyle{definition}
\newtheorem{definition}{Definition}

\newtheorem{remark}{Remark}

\numberwithin{equation}{section}

\newcommand{\qbinomial}[3]{\mbox{$
\biggl[ 
\begin{array}{c}
#1\\
 #2
\end{array}\biggr]_{
\!{#3}}$}}

\newcommand{\be}{\begin{equation}}
\newcommand{\ee}{\end{equation}}
\newcommand{\bea}{\begin{eqnarray}}
\newcommand{\eea}{\end{eqnarray}}
\newcommand{\bd}{\begin{displaymath}}
\newcommand{\ed}{\end{displaymath}}

\begin{document}

\begin{center}

\textbf{\large A General Family of $q$-Hypergeometric Polynomials \\[1mm]
and Associated Generating Functions}\\[4mm]

\textbf{H. M. Srivastava$^{1,2,3,4,\ast}$ and Sama Arjika$^{5}$}\\[2mm]

$^{1}$Department of Mathematics and Statistics, 
University of Victoria, \\
Victoria, British Columbia $V8W \; 3R4$, Canada\\

$^{2}$Department of Medical Research, 
China Medical University Hospital, \\
China Medical University,
Taichung $40402$, Taiwan, Republic of China\\

$^{3}$Department of Mathematics and Informatics, 
Azerbaijan University,\\ 
71 Jeyhun Hajibeyli Street, 
AZ$1007$ Baku, Azerbaijan\\

$^{4}$Section of Mathematics, 
International Telematic University Uninettuno,\\
I-$00186$ Rome, Italy\\[1mm]

\textbf{E-Mail: harimsri@math.uvic.ca}\\

\textbf{$^{\ast}$Corresponding Author}\\[2mm]

$^{5}$Department of Mathematics and Informatics, 
University of Agadez,\\ 
Post Office Box $199$, Agadez, Niger

\textbf{E-Mail: rjksama2008@gmail.com}

\end{center}

\vskip 4mm 

\begin{center} 
{\bf Abstract}
\end{center}
\begin{quotation}
In this paper, we introduce a general family of  
$q$-hypergeometric polynomials and investigate 
several $q$-series identities such 
as an extended generating function
and a Srivastava-Agarwal type bilinear generating function for
this family of $q$-hypergeometric polynomials. 
We give a transformational identity involving 
generating functions for the generalized 
$q$-hypergeometric polynomials
which we have introduced here. 
We also point out relevant connections
of the various $q$-results, 
which we investigate here,
with those in several related earlier works 
on this subject. We conclude this
paper by remarking that it will be a rather trivial 
and inconsequential exercise
to give the so-called $(p,q)$-variations of the $q$-results,
which we have investigated here, because the additional
parameter $p$ is obviously redundant. 
\end{quotation}   

\vskip 2mm


\vskip 1mm


\section{\bf Introduction, Definitions and Preliminaries}
\label{section1} 

In this paper, we adopt the commonly-used conventions and notations 
for the basic (or $q$-) series and the basic (or $q$-) polynomials.  
For the convenience of the reader, we first provide a summary of the 
mathematical notations, basic properties and definitions to be used  
in this paper. We refer  to the general references (see, for example,
\cite{Koekock}, \cite{Slater} and \cite{SrivastavaKarlsson})  
for the definitions and notations. 

Throughout this paper, we assume that  $|q|< 1$. For complex numbers $a$, 
the $q$-shifted factorials are defined by
\bea
(a;q)_0:=1,\quad  (a;q)_{n} =\prod_{k= 0}^{n-1} (1-aq^k),   
\quad (a;q)_{\infty}:=\prod_{k=0}^{\infty}(1-aq^{k})
\eea
and 
$$(a_1,a_2, \cdots, a_r;q)_m=(a_1;q)_m (a_2;q)_m\cdots(a_r;q)_m
\qquad (m\in\mathbb{N}_0:=\{0, 1, 2,\cdots\}=\mathbb{N}\cup\{0\}).$$  
We will frequently use the following $q$-identity \cite{Mohammed}: 
\begin{equation}
\label{used}
(aq^{-n};q)_n=(q/a;q)_n(-a)^nq^{-n-({}^n_2)}.
\end{equation}

Here, in this paper, we are mainly concerned with 
the Cauchy polynomials defined by (see \cite{Chen2003}):
\bea
\label{def}
P_n(x,y):=(x-y)(x-qy)\cdots \big(x-q^{n-1}y\big) 
=(y/x;q)_n\; x^n
\eea
together with their Srivastava-Agarwal type generating function
as giv\textcolor{red}{e}n below:  
\be
\label{Srivas}
\sum_{n=0}^\infty  
P_n (x,y)\frac{(\lambda;q)_n\,t^n}{(q;q)_n}
= {}_{2}\Phi_1\left[
\begin{array}{r}\lambda,y/x;\\
\\
0;
\end{array} 
q; \; xt\right].
\ee
For $\lambda=0,$ we get the following generating function 
(see \cite{Chen2003}):
\be
\label{gener}
\sum_{n=0}^{\infty} P_n(x,y)
\frac{t^n }{(q;q)_n} = 
\frac{(yt;q)_\infty}{(xt;q)_\infty}.
\ee

The generating function (\ref{gener}) is also 
the homogeneous version of the Cauchy identity or 
the $q$-binomial theorem  is given by 
(see \cite{GasparRahman})
\be
\label{putt}
\sum_{k=0}^{\infty} 
\frac{(a;q)_k }{(q;q)_k}z^{k}={}_{1}\Phi_0\left[
\begin{array}{r}a;\\
\\
\overline{\hspace{3mm}}\;;
\end{array} 
q;\; z\right]= 
\frac{(az;q)_\infty}{(z;q)_\infty}\qquad (|z|<1). 
\ee
Putting $a=0$, the relation (\ref{putt}) becomes 
Euler's identity  (see \cite{GasparRahman}):
\be
\label{q-expo-alpha}
\sum_{k=0}^{\infty} \frac{ z^{k}}{(q;q)_k}
=\frac{1}{(z;q)_\infty} \qquad (|z|<1)
\ee
and its inverse relation  given below 
(see \cite{GasparRahman}):
\be
\label{q-Expo-alpha}
\sum_{k=0}^{\infty}  \frac{(-1)^kq^{ ({}^k_2)
}\,z^{k}}{(q;q)_k}=(z;q)_\infty.
\ee

Saad and  Sukhi \cite{Saadsukhi} introduced the dual   
homogeneous $q$-difference operator ${\Theta}_{xy}$  
as follows (see also \cite{Liu97,SrivastaAbdlhusein,JC-SA2021C}):
\be 
\label{deffd}  
{\Theta}_{xy}\big\{f(x,y)\}:=\frac{f(q^{-1}x,y)-f( x,qy)}
{q^{-1}x-y},
\ee
which, when it acts upon 
functions of suitably-restricted variables $x$ and $y$,
yields
\begin{equation}
\Theta_{xy}^k  \{
P_n(y,x)\}=(-1)^k\frac{(q;q)_n}{(q;q)_{n-k}}\,P_{n-k}(y,x)
\quad \text{and} \quad \Theta_{xy}^k \left\{\frac{(xt;q)_\infty}
{(yt;q)_\infty}\right\}=(-t)^k\frac{(xt;q)_\infty}{(yt;q)_\infty}.\label{OP}
\end{equation}

The Hahn polynomials (see \cite{Hahn049}, 
\cite{Hahn1949} and \cite{Hahn49})
or, equivalently, the Al-Salam-Carlitz $q$-polynomials 
(see \cite{AlSalam}) are defined as follows:

\begin{align}
\phi_n^{(a)}(x|q)=\sum_{k=0}^n\begin{bmatrix}
n \\
k \\
\end{bmatrix}_q(a;q)_k\;x^k\quad \text{and} \quad
\psi_n^{(a)}(x|q)=\sum_{k=0}^n\begin{bmatrix}
n \\
k \\
\end{bmatrix}_qq^{k(k-n)}(aq^{1-k};q)_k\;x^k.\label{1.11}
\end{align}
These polynomials are usually called the 
``{\it Al-Salam-Carlitz polynomials}" in
several recent publications. Moreover, because of their
considerable role in the theories of $q$-series and 
$q$-orthogonal polynomials, many  
authors investigated various extensions of the Al-Salam-Carlitz 
polynomials (see, for example, \cite{CaoJ2013}, \cite{JC-BX-SA2020}
and \cite{HMS-Sama2020}).
 
Recently, Cao \cite[Eq. (4.7)]{CaoJ2013} introduced 
the following two families 
of generalized Al-Salam-Carlitz polynomials:
\begin{align}
\label{BE}
\phi_n^{(a,b,c)}(x,y|q)=\sum_{k=0}^n\begin{bmatrix}
n \\
k \\
\end{bmatrix}_q\frac{(a,b;q)_k}{(c;q)_k}\;x^k\; y^{n-k}
\end{align}
and

\begin{align}
\psi_n^{(a,b,c)}(x,y|q)=\sum_{k=0}^n\begin{bmatrix}
n \\
k \\
\end{bmatrix}_q\frac{(-1)^kq^{{k+1\choose2}-nk}(a,b;q)_k}
{(c;q)_k}\;x^k\;y^{n-k}, \label{BE1}
\end{align}
together with their generating functions given by 
(see \cite[Eqs.(4.10) and (4.11)]{CaoJ2013})
\begin{align}\label{9}
\sum_{n=0}^\infty \phi_n^{(a,b,c)}(x,y|q)\;\frac{t^n}{(q;q)_n}
=\frac{1}{(xt;q)_\infty}{}_2\Phi_1\left[\begin{matrix}
\begin{array}{rr}
a,b;\\
\\
c;
\end{array}
\end{matrix}\; q;yt\right]\qquad (\max\{|yt|,|xt|\}<1)
\end{align}
and

\begin{align}\label{10}
\sum_{n=0}^\infty \psi_n^{(a,b,c)}(x,y|q)
\frac{(-1)^nq^{n\choose 2}t^n}{(q;q)_n}=(xt;q)_\infty
{}_2\Phi_1\left[
\begin{matrix}
\begin{array}{rr}
a,b;\\
\\
c;
\end{array}
\end{matrix}\; q;yt\right] \qquad (|xt|<1).
\end{align}
\begin{remark}
 For $(b,c,y)=(0,0,1)$ and For $(a,b,c,y)=(1/a,0,0,1)$, (\ref{BE}) and (\ref{BE1}) reduce  to $\phi_n^{(a)}(x)$ and $\psi_n^{(a)}(x)$ defined in  (\ref{1.11}).
\end{remark}

Motivated by the above-cited work \cite{CaoJ2013},  
Cao {\it et al.} \cite{JC-BX-SA2020} introduced   
the extensions $\phi_n^{a,b,c\choose d,e}(x,y|q)$ 
and $\psi_n^{a,b,c\choose d,e}(x,y|q)$ of the
Al-Salam-Carlitz polynomials, which are defined by
\begin{align}
\phi_n^{a,b,c\choose d,e}(x,y|q)=\sum_{k=0}^n\begin{bmatrix}
n \\
k \\
\end{bmatrix}_q\frac{(a,b,c;q)_k}
{(d,e;q)_k}\; x^{n-k}\; y^k\label{23}
\end{align}
and

\begin{align}
\psi_n^{a,b,c\choose d,e}(x,y|q)=\sum_{k=0}^n\begin{bmatrix}
n \\
k \\
\end{bmatrix}_q\frac{(-1)^kq^{k(k-n)}(a,b,c;q)_k}
{(d,e;q)_k}\; x^{n-k}\; y^k,\label{24}
\end{align}
respectively. More recently, Srivastava and Arjika 
\cite{HMS-Sama2020}
introduced two families 
$\phi_n^{({\bf a},{\bf b})}(x,y|q)$ and 
$\psi_n^{({\bf a},{\bf b})}(x,y|q)$
of the generalized Al-Salam-Carlitz $q$-polynomials, 
which are defined by
\begin{equation}
\label{Carl}
\phi_n^{({\bf a}, {\bf b})}(x,y|q)
=\sum_{k=0}^n\qbinomial{n}{k}{q}\;
\frac{(a_1,a_2,\cdots,a_{r+1};q)_k}
{(b_1,b_2,\cdots,b_{r};q)_k}\; x^k\;y^{n-k}
\end{equation}
and

\begin{equation}
\label{Carl1}
\psi_n^{({\bf a},{\bf b})}(x,y|q)
=\sum_{k=0}^n \qbinomial{n}{k}{q}\;
\frac{(a_1,a_2,\cdots,a_{r+1};q)_k}
{(b_1,b_2,\cdots,b_{r};q)_k} 
q^{\binom{k+1}{2}-nk}\,x^k\;y^{n-k},
\end{equation}
respectively, where 
$${\bf a}:=(a_1,a_2,\cdots, a_{r+1})\qquad \text{and} \qquad  
{\bf b}=(b_1,b_2,\cdots, b_r).$$ On the other hand, 
Cao \cite{CaoJ2013} introduced and studied the following
family of $q$-polynomials:
\be 
\label{Varndef}
V_n^{({\bf a},{\bf c})}(x,y,z|q)
=\sum_{k=0}^n {\,n\,\atopwithdelims []\,k\,}_{q} 
\frac{(a_1,a_2,\cdots, a_r;q)_k}{(c_1,c_2,\cdots, c_u;q)_k}
\; P_{n-k}(x,y)z^k,
\ee
 where 
$${\bf a}:=(a_1,a_2,\cdots, a_{r})\qquad \text{and} \qquad  
{\bf c}=(c_1,c_2,\cdots, c_s),$$
together with their generating function given by
\be
\label{Identity}
\sum_{n=0}^\infty V_n^{({\bf a},{\bf c})}(x,y,z|q)\;\frac{t^n}{(q;q)_n}
=\frac{(yt;q)_\infty}{(xt;q)_\infty}\;{}_{r}\Phi_{\textcolor{blue}{u}}\left[
\begin{array}{rr}a_1,a_2,\cdots,a_r; \\
\\
c_1,c_2,\cdots,c_u;
 \end{array} 
\;q;zt\right].
\ee 

Our present investigation is motivated essentially by
the earlier works by Srivastava and Arjika \cite{HMS-Sama2020} 
and  Cao  \cite{CaoJ2013}.
Our aim here is to introduce and study some further extensions 
of the above-mentioned $q$-polynomials.

\begin{definition}
{\rm In terms of $q$-binomial coefficient, a family of 
generalized $q$-hypergeometric polynomials 
$\Psi_n^{({\bf a},{\bf b})}(x,y,z|q)$ are defined by
\be 
\label{ndef}
\Psi_n^{({\bf a},{\bf b})}(x,y,z|q):=(-1)^nq^{-({}^n_2)} 
\sum_{k=0}^n {\,n\,\atopwithdelims []\,k\,}_{q}
\Big[(-1)^kq^{({}^k_2)}\Big]^{1+s-r}W_k({\bf a},{\bf b})
P_{n-k}(y,x)z^k,
\ee
where, for convenience,} 
$${\bf a}=(a_1,a_2,\cdots,a_r),
\quad  
{\bf b}=(b_1,b_2,\cdots,b_s)\quad \text{and} \quad 
W_k({\bf a},{\bf b}):=\frac{(a_1,a_2,\cdots, a_r;q)_k}
{(b_1,b_2,\cdots, b_s;q)_k}.$$
\end{definition}

The above-defined $q$-polynomials 
$\Psi_n^{({\bf a},{\bf b})}(x,y,z|q)$ 
include many one-variable $q$-hypergeometric ${}_r\Phi_s$ series  
as special or limit cases. Therefore, we choose
just  to call them generalized $q$-hypergeometric polynomials.

\begin{remark}
{\rm The generalized $q$-hypergeometric polynomials 
$\Psi_n^{({\bf a},{\bf b})}(x,y,z|q)$ defined 
in $(\ref{ndef})$  are   generalized  and unified form of 
the Hahn polynomials and the Al-Salam-Carlitz polynomials.
\begin{enumerate}
\item Upon setting $(r, x, y)=(u+1, y,x)$ and ${\bf b}={\bf c}$,  
the generalized $q$-hypergeometric polynomials $(\ref{ndef})$ 
would reduce to the generalized 
Al-Salam-Carlitz $q$-polynomials 
$V_n^{({\bf a},{\bf c})}(x,y|q)$  in (\ref{Varndef}) (see \cite{CAO2016}):
\be
\Psi_n^{({\bf a},{\bf c})}(y,x,z|q)
=(-1)^nq^{-({}^n_2)}V_n^{({\bf a},{\bf c})}(x,y,z|q).
\ee
\item By choosing $r=s+1,\;x=0$ and $z=x$, the generalized  
$q$-hypergeometric polynomials $(\ref{ndef})$ reduce to the  
generalized Al-Salam-Carlitz $q$-polynomials 
$\phi_n^{({\bf a},{\bf b})}(x,y|q)$ (see \cite{HMS-Sama2020}):
\be
\Psi_n^{({\bf a},{\bf b})}(0,y,x|q)
=(-1)^nq^{-({}^n_2)}\phi_n^{({\bf a},{\bf b})}(x,y|q).
\ee
\item  Upon setting $r=s+1,\;y=0,\,z=-x$ and $x=y$, 
the $q$-hypergeometric polynomials $(\ref{ndef})$ reduce to the  
generalized Al-Salam-Carlitz $q$-polynomials 
$\psi_n^{({\bf a},{\bf b})}(x,y|q)$ (see \cite{HMS-Sama2020}):
\be
\Psi_n^{({\bf a},{\bf b})}(0,y,-x|q)
=\psi_n^{({\bf a},{\bf b})}(x,y|q).
\ee
\item  For $r=3,\,s=2,\, {\bf a}=(a,0,0),\,{\bf b} 
=(0,0)$ and $z=b$, 
the $q$-hypergeometric polynomials $(\ref{ndef})$ reduce to  
the generalized Hahn polynomials $h_n(x,y,a,b|q)$ 
(see \cite{5}):
\be
\Psi_n^{({\bf a},{\bf c})}(y,x,z|q)
=(-1)^nq^{-{n\choose2}}h_n(x,y,a,b|q).
\ee
\item  For $r=s$ and ${\bf a}={\bf b}={\bf 0}$,  
the $q$-hypergeometric polynomials 
$\Psi_n^{({\bf a},{\bf b})}(x,y,z|q)$ are the 
well known trivariate $q$-polynomials $F_n(x,y,z;q)$ 
(see \cite{Mohammed}): 
\be
\Psi_n^{({\bf 0},{\bf 0})}(x,y,z|q)=F_n(x,y,z;q):=(-1)^nq^{-({}^n_2)} 
\sum_{k=0}^n {\,n\,\atopwithdelims []\,k\,}_{q}
 (-1)^kq^{({}^k_2)} 
z^kP_{n-k}(y,x).
\ee
\item  For $r=3,\,s=2,\,{\bf a}= (a,b,c),\, {\bf b}=(d,e)$ 
and $(x,y,z)=(0,x,y)$, the generalized $q$-hypergeometric 
polynomials $\Psi_n^{({\bf a},{\bf b})}(y,x,z|q)$ reduce  
to $\phi_n^{a,b,c\choose d,e}(x,y,z|q)$ (see \cite{JC-BX-SA2020}):
\begin{equation}
\Psi_n^{({\bf a},{\bf b})}(y,x,z|q)=(-1)^nq^{-({}^n_2)} 
\phi_n^{a,b,c\choose d,e}(x,y,z|q)\textcolor{red}{.}
\end{equation}
\item  For $r=s=2,\,{\bf a}= (a,b,c),\, {\bf b}=(d,e)$  
and $(x,y,z)=(0,x,y)$, the $q$-hypergeometric polynomials 
$(\ref{ndef})$ reduce to the polynomials 
$\psi_n^{a,b,c\choose d,e}(x,y,z|q)$ (see \cite{JC-BX-SA2020}):
\begin{equation}
\Psi_n^{({\bf a},{\bf b})}(y,x,z|q)=(-1)^nq^{-({}^n_2)} 
\psi_n^{a,b,c\choose d,e}(x,y,z|q).
\end{equation}
\item If we let $r=2,\,s=1,\,{\bf a}=(a,0),\,{\bf b}=(0),\, x=0$ and $z=x$, 
the polynomials $\Psi_n^{({\bf a},{\bf b})}(x,y,z|q)$ reduce to  
the second Hahn polynomials $\phi_n^{(a)}(x,y|q)$ (see \cite{JCao12}):
\be
\Psi_n^{({\bf 0},{\bf 0})}(x,ax,y|q)=(-1)^nq^{-({}^n_2)}\phi_n^{(a)}(x,y|q).
\ee
\item If we let $r=2,\,s=1,\,({\bf a},{\bf b})=({\bf 0},{\bf 0}),\, y = ax$ 
and $z=y$, the polynomials $\Psi_n^{({\bf a},{\bf b})}(x,y,z|q)$ reduce to  
the second Hahn polynomials $\psi_n^{(a)}(x,y|q)$ (see \cite{JCao12}):
\be
\Psi_n^{({\bf 0},{\bf 0})}(x,ax,y|q)=\psi_n^{(a)}(x,y|q).
\ee
\item 
For $r=2,\,s=1,\,({\bf a},{\bf b})=({\bf 0},{\bf 0}),\, x=y=0$ and $z=x$,  
the polynomials $\Psi_n^{({\bf a},{\bf b})}(x,y,z|q)$   
reduce to the generalized Hahn polynomials $\phi_n^{(a)}(x|q)$  
(see \cite{Hahn049}, \cite{Hahn1949} and \cite{Hahn49})
or, equivalently, the Al-Salam-Carlitz $q$-polynomials 
(see \cite{AlSalam}): 
\be
\Psi_n^{({\bf 0},{\bf 0})}(0,0,x|q) =(-1)^nq^{-({}^n_2)} 
\phi_n^{(a)}(x|q). \label{ipart}
\ee
\item
Upon putting $r=s,\,({\bf a},{\bf b})=({\bf 0},{\bf 0}),\, y = ax$ and $z=1$,  
the  generalized $q$-hypergeometric polynomials 
$\Psi_n^{({\bf a},{\bf b})}(x,y,z|q)$   
reduce to the Hahn polynomials $\psi_n^{(a)}(x|q)$  
(see \cite{Hahn049}, \cite{Hahn1949} and \cite{Hahn49}) 
or, equivalently, the Al-Salam-Carlitz $q$-polynomials 
(see \cite{AlSalam}): 
\be
\Psi_n^{({\bf 0},{\bf 0})}(x,ax,1|q) = \psi_n^{(a)}(x|q). 
\label{part}
\ee
\end{enumerate}}
\end{remark}

This paper is organized as follows.
In Section \ref{Section1}, we introduce 
the following homogeneous $q$-difference   
operator: 
$${}_{r}\Phi_s\left[\begin{array}{r}a_1,a_2,\cdots,a_r; \\
\\
b_1,b_2,\cdots,b_s;
\end{array} 
\;q;-z\Theta_{xy}\right]$$ 
and apply it to investigate several $q$-series properties. 
In addition, we derive several extended generating functions 
for the generalized $q$-hypergeometric polynomials.  
In Section \ref{section3}, we first state and prove the Rogers type
formulas. The Srivastava-Agarwal type bilinear generating functions  
involving the generalized $q$-hypergeometric polynomials are 
derived in Section \ref{section4}.
In Section \ref{section5}, 
we give a transformational identity involving generating 
functions for generalized $q$-hypergeometric polynomials. Finally,
in our concluding Section \ref{section6}, we present several
remarks and observations. We also reiterate the well-documented
fact that it will be a rather trivial exercise
to give the so-called $(p,q)$-variation of the $q$-results,
which we have investigated here, because the additional
parameter $p$ is obviously redundant. 

\section{\bf Generalized $q$-Hypergeometric Polynomials}
\label{Section1}

In this section, we begin by introducing a homogeneous 
$q$-difference hypergeometric operator as follows.

\begin{definition} 
{\rm The homogeneous $q$-difference hypergeometric 
operator is defined by
\begin{equation}
\label{operator}
{}_{r}\Phi_s\left[\begin{array}{rr}a_1,a_2,\cdots,a_r;\\
\\
b_1,b_2,\cdots,b_s;
\end{array} 
\; q;-z\Theta_{xy}\right]
=\sum_{k=0}^\infty W_k({\bf a},{\bf b}) 
\frac{(-z\Theta_{xy})^k}{(q;q)_k}
\left[(-1)^kq^{({}^k_2)}\right]^{1+s-r}.
\end{equation}}
\end{definition}

We now derive the $q$-series identities (\ref{conds}), 
(\ref{identity}) and (\ref{bella}) below,  
which will be used later in order to derive 
the extended generating functions, 
the Rogers type formulas and the Srivastava-Agarwal type 
bilinear generating functions involving the generalized 
$q$-hypergeometric polynomials.

\begin{lemma}
\label{MAL} 
Each of the following operational formulas holds true
for the homogeneous $q$-difference operator defined by
$\eqref{operator}$$:$
\begin{align}
\label{conds}
{}_{r}\Phi_s\left[\begin{array}{r}a_1,a_2,\cdots,a_r;\\
\\
b_1,b_2,\cdots,b_s;
\end{array} 
\;q;-z\Theta_{xy}\right]\left\{(-1)^nq^{-({}^n_2)} P_n(y,x)\right\}
=\Psi_n^{({\bf a},{\bf b})}(x,y,z|q),
\end{align}

\begin{align}
\label{identity}
{}_{r}\Phi_s\left[\begin{array}{r}a_1,a_2,\cdots,a_r;\\
\\
b_1,b_2,\cdots,b_s;
\end{array} 
\;q;-z\Theta_{xy}\right]\left\{\frac{(xt,q)_\infty}
{(yt;q)_\infty}\right\}
=\frac{(xt;q)_\infty}{(yt;q)_\infty}\;{}_{r}\Phi_s\left[
\begin{array}{rr}a_1,a_2,\cdots,a_r;\\
\\
b_1,b_2,\cdots,b_s;
\end{array} 
q;zt\right]
\end{align}
and

\begin{align}
\label{bella}
&{}_{r}\Phi_s\left[\begin{array}{r}a_1,a_2,\cdots,a_r;\\
\\
b_1,b_2,\cdots,b_s;
\end{array} 
\;q;-z\Theta_{xy}\right]\left\{\frac{P_k(y,x)\,(xt;q)_\infty}
{(xt;q)_k(yt;q)_\infty}\right\}\cr
&\qquad\qquad=t^{-k}\frac{ (xt;q)_\infty}{ (yt;q)_\infty}
\sum_{j=0}^k\frac{ (q^{-k},yt;q)_j\,q^j}{(xt,q;q)_j} 
\;{}_{r}\Phi_s\left[\begin{array}{rr}a_1,a_2,\cdots,a_r;\\ 
\\
b_1,b_2,\cdots,b_s;
\end{array}
\; q; ztq^j\right]\qquad (|yt|<1).
\end{align} 
\end{lemma}

\begin{proof}
Firstly, by applying (\ref{operator}) and (\ref{OP}), 
we get (\ref{conds}). 
Secondly, in light of (\ref{OP}),  
we get the desired identity (\ref{identity}).  
Thirdly, by making use of the
$q$-Chu-Vandermonde  formula \cite[Eq. (II.6)]{GasparRahman}
\be
\label{qchuv}
{}_{2}\Phi_1\left[\begin{array}{rr} q^{-n},a;\\
\\
c;
\end{array} 
\; q; q \right]=\frac{(c/a;q)_n}{(c;q)_n}\;a^n,
\ee
we find that
\begin{equation}
\label{equation}
\frac{P_k(y,x)\,(xt;q)_\infty}{(xt;q)_k(yt;q)_\infty}
=t^{-k}\; \frac{(xt;q)_\infty}{(yt;q)_\infty} 
{}_{2}\Phi_1\left[\begin{array}{rr} q^{-k}, yt;\\
\\
xt;
\end{array} 
\; q; q \right]=t^{-k}\sum_{j=0}^\infty\frac{(q^{-k};q)_j\,q^j}
{(q;q)_j}\; \frac{(xtq^j;q)_\infty}{(ytq^j;q)_\infty}.
\end{equation}
Therefore, by using \eqref{equation} and \eqref{identity} 
successively, we obtain
\begin{align*}
& {}_{r}\Phi_s\left[\begin{array}{rr}a_1,a_2,\cdots,a_r;\\
\\
b_1,b_2,\cdots,b_s;
\end{array} 
\;q;-z\Theta_{xy}\right]\left\{\frac{P_k(y,x)\,(xt;q)_\infty}
{(xt;q)_k(yt;q)_\infty}\right\}\notag \\
&\qquad\qquad= t^{-k}\sum_{j=0}^k\frac{ (q^{-k};q)_j\,q^j}{(q;q)_j}  
{}_{r}\Phi_s\left[
\begin{array}{rr}a_1,a_2,\cdots,a_r;\\
\\
b_1,b_2,\cdots,b_s;
\end{array} 
\;q;-z\Theta_{xy}\right]\left\{\frac{(xtq^j;q)_\infty}
{(ytq^j;q)_\infty} 
\right\} \notag \\
&\qquad\qquad=  t^{-k}\sum_{j=0}^k\frac{ (q^{-k};q)_j\,q^j}{(q;q)_j} 
\frac{(xtq^j;q)_\infty}{(ytq^j;q)_\infty}\;{}_{r}\Phi_s\left[
\begin{array}{r}a_1,a_2,\cdots,a_r;\\
\\
b_1,b_2,\cdots,b_s;
\end{array} 
\; q; ztq^j\right], 
\end{align*}
which is asserted by Lemma \ref{MAL}.
\end{proof}

We now derive an extended generating function for 
the generalized $q$-hypergeometric polynomials  
$\Psi_{n}^{({\bf a},{\bf b})}(x,y|q)$ by using the 
operator representation (\ref{conds}).

\begin{theorem}[Extended generating function] 
For $k\in\mathbb{N}$ and $\max\{|yt|,|zt|\}<1,$ 
it is asserted that 
\label{exgend}
\begin{align}
\label{exgen}
&\sum_{n=0}^\infty\Psi_{n+k}^{({\bf a},{\bf b})}(x,y,z|q) 
(-1)^{n+k}q^{({}^{n+k}_{\;\;\;2})}\;
\frac{t^n }{(q;q)_n}\cr
&\qquad\qquad=\frac{(xt;q)_\infty}{t^k(yt;q)_\infty}
\sum_{j=0}^k\frac{ (q^{-k},yt;q)_j\,q^j}{(q,xt;q)_j} 
\;{}_{r}\Phi_s\left[\begin{array}{rr}a_1,a_2,\cdots,a_r;\\
\\
b_1,b_2,\cdots,b_s;
\end{array} 
\; q; ztq^j\right] \qquad (|yt|<1).
\end{align}
\end{theorem}

\begin{proof} 
In light of (\ref{conds}) and \eqref{gener}, we have 
\begin{align*}
&\sum_{n=0}^\infty\Psi_n^{({\bf a},{\bf b})}(x,y,z|q)
(-1)^{n+k}q^{ ({}^{n+k}_{\;\;\;2})}\frac{ t^n}{(q;q)_n} \\
&\qquad\qquad= 
\sum_{n=0}^\infty {}_{r}\Phi_s\left[
\begin{array}{rr}a_1,a_2,\cdots,a_r;\\
\\
b_1,b_2,\cdots,b_s;
\end{array} 
\;q;-z\Theta_{xy}\right]\left\{ (-1)^{n+k}q^{ ({}^{n+k}_{\;\;\;2})} 
P_{n+k}(y,x)\right\}\frac{ (-1)^{n+k}q^{ ({}^{n+k}_{\;\;\;2})}\,t^n}{(q;q)_n}\cr
&\qquad\qquad={}_{r}\Phi_s\left[\begin{array}{rr}a_1,a_2,\cdots,a_r;\\
\\
b_1,b_2,\cdots,b_s;
\end{array} 
\;q;-z\Theta_{xy}\right]\left\{P_k(y,x)\sum_{n=0}^\infty  
P_n(y,xq^k)\frac{t^n}{(q;q)_n}\right\} \\
&\qquad\qquad= {}_{r}\Phi_s\left[
\begin{array}{rr}a_1,a_2,\cdots,a_r;\\
\\
b_1,b_2,\cdots,b_s;
\end{array} 
\;q;-z\Theta_{xy}\right]\left\{\frac{P_k(y,x)\,(xt;q)_\infty}
{(xt;q)_k(yt;q)_\infty}\right\}.
\end{align*}
The proof of Theorem \ref{exgend} is completed by using (\ref{bella}).
\end{proof}

\begin{remark}
Setting $k=0$ in Theorem \ref{exgend}, we get the 
following generating function for the generalized 
$q$-hypergeometric polynomials:
\begin{equation}
\label{gen}
\sum_{n=0}^\infty  \Psi_n^{({\bf a},{\bf b})}(x,y,z|q) 
(-1)^nq^{ ({}^n_2)}\frac{\,t^n}{(q;q)_n}
=\frac{(xt;q)_\infty}{(yt;q)_\infty}\;{}_{r}\Phi_s\left[
\begin{array}{rr}a_1,a_2,\cdots,a_r;\\
\\
b_1,b_2,\cdots,b_s;
\end{array} 
\; q; zt\right].
\end{equation}
\end{remark}

\section{\bf The Rogers Formula}
\label{section3}

In this section, we use the assertion (\ref{bella}) 
of Lemma \ref{MAL} in order to derive several 
$q$-identities such as the Rogers type formula for 
the generalized $q$-hypergeometric polynomials   
$\Psi_{n}^{({\bf a},{\bf b})}(x,y|q)$.

\begin{theorem}
[The Rogers formula for $\Psi_{n}^{({\bf a},{\bf b})}(x,y,z|q)$]
\label{dsexgend} 
For $\max\left\{\left|\frac{t}{\omega}\right|,|y\omega|\right\}<1,$
the following Rogers type  formula holds$:$
\begin{align} 
\label{zdd}
&\sum_{n=0}^\infty \sum_{k=0}^\infty \Psi_{n+k}^{({\bf a},{\bf b})}
(x,y,z|q)(-1)^{n+k}q^{ ({}^{n+k}_{\,\,\,2})}\,\frac{t^n}{(q;q)_n}\;
\frac{\omega^k}{(q;q)_k}\cr
&\qquad\qquad=\frac{(x\omega;q)_\infty}{(t/\omega,y\omega;q)_\infty} 
\sum_{k=0}^\infty\frac{ (y \omega;q)_k\,q^k}{(q\omega/t,x\omega,q;q)_k} 
\;{}_{r}\Phi_s\left[\begin{array}{rr}a_1,a_2,\cdots,a_r;\\ 
\\
b_1,b_2,\cdots,b_s;
\end{array} 
\; q; z\omega q^k\right].
\end{align} 
\end{theorem}

\begin{proof}
In light of (\ref{conds}), we have 
\begin{align} 
&\sum_{n=0}^\infty \sum_{k=0}^\infty 
\Psi_{n+k}^{({\bf a},{\bf b})}(x,y,z|q) 
(-1)^{n+k}q^{ ({}^{n+k}_{\,\,\,2})} \frac{t^n}{(q;q)_n}\;
\frac{ \omega^k}{(q;q)_k}\cr
&\qquad\qquad= 
\sum_{n=0}^\infty\sum_{k=0}^\infty{}_{r}\Phi_s\left[
\begin{array}{c}a_1,a_2,\cdots,a_r;\\ 
\\
b_1,b_2,\cdots,b_s;
\end{array} 
\; q;- z\Theta_{xy}\right]\left\{P_{n+k}(y,x)\right\} 
\frac{t^n}{(q;q)_n}\frac{\omega^k}{(q;q)_k}\cr
&\qquad\qquad={}_{r}\Phi_s\left[
\begin{array}{r}a_1,a_2,\cdots,a_r;\\ 
\\
b_1,b_2,\cdots,b_s;
\end{array} 
\; q;- z\Theta_{xy}\right]\left\{\sum_{n=0}^\infty  
P_n(y,x)\frac{ t^n}{(q;q)_n}\sum_{k=0}^\infty  
P_k(y,q^nx)\frac{\omega^k}{(q;q)_k}\right\}\cr
&\qquad\qquad={}_{r}\Phi_s\left[
\begin{array}{r}a_1,a_2,\cdots,a_r;\\ 
\\
\; b_1,b_2,\cdots,b_s;
\end{array} 
\; q;- z\Theta_{xy}\right]\left\{\sum_{n=0}^\infty  
P_n(y,x)\frac{t^n}{(q;q)_n}\; \frac{(x\omega q^n;q)_\infty}
{(y\omega;q)_\infty} \right\}\cr
&\qquad\qquad=\sum_{n=0}^\infty \frac{t^n}{(q;q)_n}{}_{r}
\Phi_s\left[\begin{array}{r}a_1,a_2,\cdots,a_r;\\ 
\\
b_1,b_2,\cdots,b_s;
\end{array} 
\; q;- z\Theta_{xy}\right]\left\{\frac{P_n(y,x)\,(x\omega;q)_\infty}
{(x\omega;q)_n(y\omega;q)_\infty}\right\}.
\end{align}
 
Now, by using (\ref{bella}) as well as \eqref{used}, this last relation takes 
the following form:
\begin{align} 
&\sum_{n=0}^\infty \sum_{k=0}^\infty \Psi_{n+k}^{({\bf a},{\bf b})}
(x,y,z|q) (-1)^{n+k}\;q^{({}^{n+k}_{\,\,\,\,2})}\frac{t^n}{(q;q)_n}\;
\frac{ \omega^k}{(q;q)_k}\notag \\
&\qquad\qquad=\frac{(x\omega;q)_\infty}{(y\omega;q)_\infty}
\sum_{n=0}^\infty \frac{(t/\omega)^n}{(q;q)_n}  
\sum_{k=0}^n\frac{ (q^{-n},y\omega;q)_k\,q^k}{(x\omega,q;q)_k} 
\;{}_{r}\Phi_s\left[
\begin{array}{r}a_1,a_2,\cdots,a_r;\\
\\
b_1,b_2,\cdots,b_s;
\end{array}
\; q; z\omega q^k\right]\notag \\
&\qquad\qquad=\frac{(x\omega;q)_\infty}{(y\omega;q)_\infty} 
\sum_{k=0}^\infty\frac{(y\omega;q)_k\,q^k}{(x\omega,q;q)_k} 
\;{}_{r}\Phi_s\left[
\begin{array}{r}a_1,a_2,\cdots,a_r;\\
\\
b_1,b_2,\cdots,b_s;
\end{array} 
\; q; z\omega q^k\right]\sum_{n=k}^\infty 
\frac{(t/\omega)^n(-1)^kq^{({}^k_2)-nk}}{(q;q)_{n-k}}\notag \\
&\qquad\qquad=\frac{(x\omega;q)_\infty}{(y\omega;q)_\infty} 
\sum_{k=0}^\infty\frac{(y\omega;q)_k\,(-t/\omega)^k q^{-({}^k_2)}}
{(x\omega,q;q)_k} \;{}_{r}\Phi_s\left[
\begin{array}{r}a_1,a_2,\cdots,a_r;\\
\\
b_1,b_2,\cdots,b_s;
\end{array} 
\; q; z\omega q^k\right]\sum_{n=0}^\infty 
\frac{(t/\omega)^n\;q^{-nk}}{(q;q)_{n}} \notag \\
&\qquad\qquad=\frac{(x\omega;q)_\infty}{(y\omega;q)_\infty} 
\sum_{k=0}^\infty\frac{ (y\omega;q)_k\,(-t/\omega)^k \;
q^{-({}^k_2)}}{(x\omega,q;q)_k} \;{}_{r}\Phi_s\left[
\begin{array}{rr}a_1,a_2,\cdots,a_r;\\
\\
b_1,b_2,\cdots,b_;
\end{array} 
\; q; z\omega q^k\right] \frac{1}{(tq^{-k}/\omega;q)_\infty}\notag \\
&\qquad\qquad=\frac{(x\omega;q)_\infty}{(t/\omega,y\omega;q)_\infty} 
\sum_{k=0}^\infty\frac{(y\omega;q)_k\,(-t/\omega)^k q^{-({}^k_2)}}
{(x\omega,tq^{-k}/\omega,q;q)_k} \;{}_{r}\Phi_s\left[
\begin{array}{rr}a_1,a_2,\cdots,a_r;\\
\\
b_1,b_2,\cdots,b_s;
\end{array}
\; q; z\omega q^k\right]\notag \\
&\qquad\qquad =\frac{(x\omega;q)_\infty}{(t/\omega,y\omega;q)_\infty} 
\sum_{k=0}^\infty\frac{ (y\omega;q)_k\,q^k}{(q\omega/t,x \omega,q;q)_k} 
\;{}_{r}\Phi_s\left[
\begin{array}{rr}a_1,a_2,\cdots,a_r;\\ 
\\
b_1,b_2,\cdots,b_s;
\end{array} 
\; q; z\omega q^k\right],
\end{align} 
which evidently completes the proof of Theorem \ref{dsexgend}.
\end{proof}
 
\section{\bf The Srivastava-Agarwal Type Bilinear Generating 
Functions for the Generalized $q$-Hypergeometric Polynomials 
$\Psi_n^{({\bf a},{\bf b})}(x,y,z|q)$}
\label{section4}

In this section, by applying the 
following homogeneous $q$-difference hypergeometric operator \cite{AS-SMK2019}: 
$${}_{r}\Phi_s\left[
\begin{array}{r}a_1,a_2,\cdots,a_r;\\
\\
b_1,b_2,\cdots,b_s;
\end{array} 
q;-z\Theta_{xy}\right],$$
we derive the Srivastava-Agarwal type generating 
functions for the generalized $q$-hypergeometric polynomials
$\Psi_n^{({\bf a},{\bf b})}(x,y,z|q)$ defined by \eqref{ndef}. 
We also deduce a bilinear generating function for the 
Al-Salam-Carlitz polynomials $\psi_n^{(\alpha)}(x|q)$ 
as an application of the Srivastava-Agarwal type generating functions.

\begin{lemma}{\rm (see \cite[Eq. (3.20)]{SrivastavaAgarwal} 
and \cite[Eq. (5.4)]{Cao2010A})} 
\label{LEMMA41}
Each of the following generating relations holds true$:$
\begin{align} 
\label{21sums}
\sum_{n=0}^\infty \phi_n^{(\alpha)}(x|q) (\lambda;q)_n\frac{  t^n}{(q;q)_n} 
=\frac{(\lambda t; q)_\infty }{(t;q)_\infty}   {}_2\Phi_1\left[
\begin{array}{rr} \lambda, \alpha;\\
\\
\lambda  t; 
\end{array}\,q; xt   
\right] \qquad \big(\max\{|t|, |xt|\}<1\big)
\end{align}
and

\begin{align}
\label{c1sums}
\sum_{n=0}^\infty \psi_n^{(\alpha)}(x|q) (1/\lambda;q)_n\;
\frac{(\lambda tq)^n}{(q;q)_n} 
=\frac{(xtq; q)_\infty }{(\lambda xtq;q)_\infty}{}_2\Phi_1\left[
\begin{array}{rr}1/ \lambda,1/(\alpha x);\\
\\
1/(\lambda  xt); 
\end{array}\;q; \alpha q  
\right] 
\end{align}
$$\big(\max\{|\lambda x tq|,|\alpha q|\}<1\big).$$
\end{lemma}
  
For more information about the Srivastava-Agarwal type 
generating functions for the Al-Salam-Carlitz polynomials, 
 one may refer to \cite{SrivastavaAgarwal} and \cite{JCao12}.

We now state and prove the Srivastava-Agarwal type bilinear 
generating functions asserted by Theorem \ref{AUTRE} below. 

\begin{theorem}
\label{AUTRE} 
Suppose that $\max\{|\alpha q|,|vxtq|,|xztq|\}<1$. Then
\begin{align}
\label{AUTRE0}
&\sum_{n=0}^\infty \psi_n^{(\alpha)}(x|q)   
\Psi_n^{({\bf a}, {\bf b})}(u,v,z|q)  
\frac{(-1)^nq^{({}^{n+1}_{\,\,\,\,\,2})} \,t^n}{(q;q)_n}\notag \\
&\qquad=\frac{(q/x,uxtq; q)_\infty }{(\alpha q,v  xtq;q)_\infty}
\sum_{n=0}^\infty \frac{(-1)^nq^{({}^n_2)} (1/(\alpha x),1/(uxt);q)_n}
{(q/x,1/(v xt),q; q)_n }\left(\frac{\alpha uq}{v}\right)^n \notag \\
&\qquad \qquad \qquad \cdot {}_{r}\Phi_s\left[
\begin{array}{r}a_1,a_2,\cdots,a_r;\\
\\
b_1,b_2,\cdots,b_s;
\end{array} 
\; q;xztq^{1-n}\right].
\end{align}
\end{theorem}

\begin{proof}
In order to prove Theorem \ref{AUTRE}, we need the
$q$-Chu-Vandermonde summation theorem given by  
(see \cite[Eq. (II.7)]{GasparRahman})
\be
\label{male}
{}_2\Phi_1\left[
\begin{array}{rr}q^{-n},a;\\
\\
c; 
\end{array}\; q; \frac{cq^n}{a}
\right]=\frac{(c/a;q)_n}{(c;q)_n} 
\qquad (n\in \mathbb{N}_0).
\ee

Upon letting $(\lambda,t)=(v/u,tu)$ in the equation (\ref{c1sums}),  
we obtain
\bea
\label{llms}
&&\sum_{n=0}^\infty \psi_n^{(\alpha)}(x|q)p_n(v,u)\frac{(qt)^n}{(q;q)_n}\cr
&&\qquad =\frac{(uxtq; q)_\infty}{(vxtq;q)_\infty}{}_2\Phi_1\left[
\begin{array}{rr}1/ (\alpha x), u/v;\\
\\
1/(vxt); 
\end{array}\;q; \alpha q  
\right] 
\cr
&&\qquad=\frac{(uxtq; q)_\infty }{(v  xtq;q)_\infty} 
\sum_{k=0}^\infty\frac{(1/(\alpha x);q)_k (\alpha q)^k}{(q;q)_k}\;
\frac{(u/v;q)_k   }{ (1/(v  xt);q)_k }\mbox{ by } (\ref{male})\cr
&&\qquad =\frac{(uxtq; q)_\infty}{(vxtq;q)_\infty} 
\sum_{k=0}^\infty\frac{(1/(\alpha x);q)_k 
(\alpha q)^k}{(q;q)_k}  
{}_2\Phi_1\left[
\begin{array}{rr}q^{-k}, 1/ (uxt);\\
\\
1/(v  xt); 
\end{array}\;q; \frac{uq^k}{v} 
\right] \cr
&&\qquad=\sum_{k=0}^\infty\frac{(1/(\alpha x);q)_k 
(\alpha q)^k}{(q;q)_k}\sum_{n=0}^\infty
\frac{(q^{-k};q)_n  \,q^{nk}}{(q;q)_n}
\; \frac{(1/(uxt); q)_n(uxtq; q)_\infty }
{(1/(vxt);q)_n(vxtq;q)_\infty}\left(\frac{u}{v}\right)^n \cr
&&\qquad=\sum_{k=0}^\infty\frac{(1/(\alpha x);q)_k (\alpha q)^k}
{(q;q)_k}\sum_{n=0}^\infty\frac{(q^{-k};q)_n \, q^{nk}}{(q;q)_n}
\; \frac{(uxtq^{1-n}; q)_\infty}{(vxtq^{1-n};q)_\infty}. 
\eea

By appealing appropriately to \eqref{llms}, 
the left-hand side of \eqref{AUTRE0} becomes  
\bea
&&\sum_{n=0}^\infty\psi_n^{(\alpha)}(x|q){}_{r}\Phi_s\left[
\begin{array}{r}a_1,a_2,\cdots,a_r;\\
\\
b_1,b_2,\cdots,b_s;
\end{array} 
\; q;-z\Theta_{uv}\right]\left\{p_{n}(v,u)\right\}  
\frac{(qt)^{n}}{(q;q)_n}
\cr
&&\qquad = {}_{r}\Phi_s\left[
\begin{array}{rr}a_1,a_2,\cdots,a_r;\\
\\
b_1,b_2,\cdots,b_s;
\end{array} 
\; q;-z\Theta_{uv}\right] 
\left\{\sum_{n=0}^\infty \psi_n^{(\alpha)}(x|q) p_{n}(v,u) 
\frac{(qt)^{n}}{(q;q)_n}\right\}
\cr
&&\qquad = {}_{r}\Phi_s\left[
\begin{array}{r}a_1,a_2,\cdots,a_r;\\
\\
b_1,b_2,\cdots,b_s;
\end{array} 
\; q;-z\Theta_{uv}\right] \left\{\sum_{k=0}^\infty
\frac{(1/(\alpha x);q)_k (\alpha q)^k}{(q;q)_k}
\sum_{n=0}^\infty\frac{(q^{-k};q)_n \, q^{nk}}{(q;q)_n}
\;\frac{(uxtq^{1-n}; q)_\infty }{(vxtq^{1-n};q)_\infty}\right\}
\cr
&&\qquad=\sum_{k=0}^\infty\frac{(1/(\alpha x);q)_k (\alpha q)^k}
{(q;q)_k}\sum_{n=0}^\infty\frac{(q^{-k};q)_n \, q^{nk}}{(q;q)_n} 
{}_{r}\Phi_s\left[
\begin{array}{rr}a_1,a_2,\cdots,a_r;\\
\\
b_1,b_2,\cdots,b_s;
\end{array} 
\; q;-z\Theta_{uv}\right] \left\{\frac{(uxtq^{1-n}; q)_\infty}
{(vxtq^{1-n};q)_\infty}\right\}
\cr
&&\qquad=\sum_{k=0}^\infty\frac{(1/(\alpha x);q)_k (\alpha q)^k}{(q;q)_k}
\sum_{n=0}^\infty\frac{(q^{-k};q)_n \, q^{nk}}{(q;q)_n}   
\; \frac{(uxtq^{1-n}; q)_\infty}{(vxtq^{1-n};q)_\infty}   
{}_{r}\Phi_s\left[\begin{array}{rr}a_1,a_2,\cdots,a_r;\\
\\
b_1,b_2,\cdots,b_s;
\end{array} 
\; q;xztq^{1-n}\right]
\cr
&&\qquad=\frac{(uxtq; q)_\infty}{(vxtq;q)_\infty} 
\sum_{k=0}^\infty\frac{(1/(\alpha x);q)_k (\alpha q)^k}{(q;q)_k}
\sum_{n=0}^\infty \frac{(q^{-k},uxtq^{1-n},q)_n \, q^{nk}}
{(vxtq^{1-n},q; q)_n } {}_{r}\Phi_s\left[
\begin{array}{rr}a_1,a_2,\cdots,a_r;\\
\\
b_1,b_2,\cdots,b_s;
\end{array} 
\; q;xztq^{1-n}\right]
  \cr
&&\qquad=\frac{(uxtq; q)_\infty}{(vxtq;q)_\infty}
\sum_{n=0}^\infty\frac{(-1)^nq^{({}^n_2)}(1/(uxt);q)_n}
{(1/(v  xt),q; q)_n}\left(\frac{u }{v}\right)^n  
{}_{r}\Phi_s\left[\begin{array}{rr}a_1,a_2,\cdots,a_r;\\
\\
b_1,b_2,\cdots,b_s;
\end{array} 
\; q;xztq^{1-n}\right]\cr
&&\qquad\qquad \qquad \qquad \cdot \sum_{k=n}^\infty
\frac{(1/(\alpha x);q)_k (\alpha q)^k}{(q;q)_{k-n}} 
\cr
&&\qquad=\frac{(uxtq; q)_\infty}{(vxtq;q)_\infty}
\sum_{n=0}^\infty   \frac{(-1)^nq^{({}^n_2)} (1/(\alpha x),1/(uxt);q)_n }
{(1/(vxt), q; q)_n}\left(\frac{\alpha u q }{v}\right)^n 
{}_{r}\Phi_s\left[\begin{array}{r}a_1,a_2,\cdots,a_r;\\
\\
b_1,b_2,\cdots,b_s;
\end{array} 
\; q;xztq^{1-n}\right]\cr
&&\qquad\qquad \qquad \qquad \cdot  
\sum_{k=0}^\infty\frac{(q^n/(\alpha x);q)_k (\alpha q)^k}{(q;q)_{k}}\cr
&&\qquad=\frac{(uxtq; q)_\infty}{(vxtq;q)_\infty}
\sum_{n=0}^\infty\frac{(-1)^nq^{({}^n_2)} (1/(\alpha x),1/(uxt);q)_n}
{(1/(vxt),q; q)_n }\left(\frac{\alpha uq}{v}\right)^n  
{}_{r}\Phi_s\left[
\begin{array}{r}a_1,a_2,\cdots,a_r;\\
\\
b_1,b_2,\cdots,b_s;
\end{array} 
\; q;xztq^{1-n}\right]\cr
&&\qquad\qquad\qquad \qquad \cdot \frac{(q^{1+n}/x;q)_\infty}
{(\alpha q;q)_{\infty}}
\cr
&&\qquad=\frac{(q/x,uxtq; q)_\infty}{(\alpha q,vxtq;q)_\infty}
\sum_{n=0}^\infty\frac{(-1)^nq^{({}^n_2)} (1/(\alpha x),1/(uxt);q)_n}
{(q/x,1/(vxt),q; q)_n}\left(\frac{\alpha uq }{v}\right)^n 
{}_{r}\Phi_s\left[
\begin{array}{rr}a_1,a_2,\cdots,a_r;\\
\\
b_1,b_2,\cdots,b_s;
\end{array} 
\; q;xztq^{1-n}\right], \nonumber
\eea
which is precisely the right-hand side of the   
assertion (\ref{AUTRE0}) of Theorem \ref{AUTRE}.   
The proof of Theorem \ref{AUTRE} is thus completed.
\end{proof}

\begin{remark}
{\rm In view of the special case $(\ref{part})$ of 
$\; \Psi_n^{({\bf a}, {\bf b})}(u,v,z|q),\;$
we deduce the bilinear generating function for 
the Hahn polynomials $\psi_n^{(\alpha)}(x|q)$ as
asserted by Corollary \ref{COROL} below.}
\end{remark}

\begin{corollary}
\label{COROL}
Let $\max\{|\alpha q|,|axytq|\}<1.$ Then
\begin{align}
\label{AUTREC1}
&\sum_{n=0}^\infty \psi_n^{(\alpha)}(x|q)   
\psi_n^{(a)}(y|q)\frac{(-1)^nq^{({}^{n+1}_{\,\,\,\,\,2})} \,t^n}
{(q;q)_n}\notag \\
&\qquad \qquad  =\frac{(q/x,xytq,xtq; q)_\infty}
{(\alpha q,axytq;q)_\infty}  
{}_{3}\Phi_2\left[
\begin{array}{rr}1/(\alpha x),1/(xyt),1/(xt);\\
\\
q/x, 1/(axyt);
\end{array} 
\; q;\frac{\alpha xtq}{a}\right].
\end{align}
\end{corollary}

\begin{remark}
{\rm For $r=s,\, {\bf a}={\bf b}={\bf 0},\, x=y,\,y = ay$ and $z=1,$ 
the assertion $(\ref{AUTRE0})$ reduces to Corollary \ref{COROL}.}  
\end{remark}

\section{\bf A Transformational Identity Involving Generating Functions 
for the Generalized $q$-Hypergeometric Polynomials}
\label{section5}

In this section, we derive the following transformational 
identity involving generating functions for the 
generalized $q$-hypergeometric polynomials. Once again, in our derivation, 
we apply the homogeneous $q$-difference operator (\ref{operator}).
  
\begin{theorem}
\label{asss}
Let $A(n)$ and $B(n)$ satisfy the following relationship$:$
\be
\label{sam}
\sum_{n=0}^\infty A(n) P_n(v,u)
=\sum_{n=0}^\infty B(n) \frac{(xutq^{1-n};q)_\infty}
{(xvtq^{1-n};q)_\infty}.
\ee
Then
\begin{align}
\label{samm}
\sum_{n=0}^\infty (-1)^nq^{({}^{n}_{2})} A(n) 
\Psi_n^{({\bf a},{\bf b})}(u,v,z|q)
=\sum_{n=0}^\infty B(n) \frac{(xutq^{1-n};q)_\infty}
{(xvtq^{1-n};q)_\infty} {}_{r}\Phi_s\left[
\begin{array}{rr}a_1,a_2,\cdots,a_r;\\
\\
b_1,b_2,\cdots,b_s;
\end{array} 
\; q;xztq^{1-n}\right],
\end{align}
provided that each of the series in $(\ref{sam})$ 
and  $(\ref{samm})$ is absolutely convergent.
\end{theorem}

\begin{proof} 
If we denote by $f(x,u,v,z)$ the right-hand side of (\ref{samm}), 
we find by using \eqref{sam} and \eqref{conds} that 
\begin{align*}
f(x,u,v,z)
&=\sum_{n=0}^\infty B(n)\;  {}_{r}\Phi_s\left[
\begin{array}{rr}a_1,a_2,\cdots,a_r;\\ \\
b_1,b_2,\cdots,b_s;
\end{array} 
\; q;-z\Theta_{uv}\right]\left\{\frac{(xutq^{1-n};q)_\infty}
{(xvtq^{1-n};q)_\infty}\right\}\cr
&={}_{r}\Phi_s\left[\begin{array}{rr}a_1,a_2,\cdots,a_r;\\
\\
b_1,b_2,\cdots,b_s;
\end{array} 
\; q;-z\Theta_{uv}\right]\left\{\sum_{n=0}^\infty B(n) 
\frac{(xutq^{1-n};q)_\infty}{(xvtq^{1-n};q)_\infty}\right\}\cr
&={}_{r}\Phi_s\left[
\begin{array}{r}a_1,a_2,\cdots,a_r; 
 \\\\
b_1,b_2,\cdots,b_s;
\end{array} 
\; q;-z\Theta_{uv}\right]\left\{\sum_{n=0}^\infty 
A(n)P_n(v,u)\right\}\cr
&= \sum_{n=0}^\infty A(n){}_{r}\Phi_s\left[
\begin{array}{rr}a_1,a_2,\cdots,a_r;\\
\\
b_1,b_2,\cdots,b_s;
\end{array} 
\;q;-z\Theta_{uv}\right]\left\{P_n(v,u)\right\}\cr
&=\sum_{n=0}^\infty (-1)^nq^{({}^{n}_{2})} A(n)  
\Psi_n^{({\bf a}, {\bf b})}(u,v,z|q),
\end{align*}
which is precisely the left-hand side of (\ref{samm}). 
Our demonstration of Theorem \ref{asss} is thus completed.
\end{proof}

\begin{remark}
{\rm  Upon setting  
$u=1, v=\lambda$,  
\begin{equation}\label{AnBn}
A(n)=\psi_n^{(\alpha)}(x|q)\frac{(qt)^n}{(q;q)_n}
\qquad \text{and}
\qquad 
B(n)=\sum_{k=0}^\infty\frac{(1/(\alpha x);q)_k (\alpha q)^k}
{(q;q)_k}\; 
\frac{(q^{-k};q)_n \, q^{nk}}{(q;q)_n}
\end{equation}
 in (\ref{sam}), we get $(\ref{c1sums})$.
  Moreover, 
Upon specializing $A(n)$ and $B(n)$ in $(\ref{samm})$ 
as in the last equation \eqref{AnBn}, we get the 
assertion $(\ref{AUTRE0})$ of Theorem \ref{AUTRE}.}
\end{remark} 

\section{\bf Concluding Remarks and Observations}
\label{section6}

In our present investigation, 
we have introduced a general family of  
$q$-hypergeometric polynomials and we have derived
several $q$-series identities such 
as an extended generating function
and Srivastava-Agarwal type bilinear generating functions for
this family of $q$-hypergeometric polynomials. 
We have presented a transformational identity involving 
generating functions for the generalized 
$q$-hypergeometric polynomials
which we have introduced here. 
We have also pointed out relevant connections
of the various $q$-results, 
which we have investigated in this paper,
with those in several related earlier 
works on this subject. 

We conclude this paper by remarking that, 
in the recently-published 
survey-cum-expository review article by Srivastava \cite{[e]}, 
the so-called $(p,q)$-calculus 
was exposed to be a rather trivial and inconsequential variation of 
the classical $q$-calculus, the additional parameter $p$ being redundant 
or superfluous (see, for details, \cite[p. 340]{[e]}). This
observation by Srivastava \cite{[e]} will indeed apply to any
attempt to produce the rather straightforward $(p,q)$-variations of the results
which we have presented in this paper.\\[2mm]

\noindent
{\bf Conflicts of Interest:} The authors declare that they 
have no conflicts of interest.

\end{document}